\documentclass[twoside,12pt,a4paper]{article}
{\markboth{\sl  Hichem H}{\sl weakly converges to }\thispagestyle{plain}
\usepackage{amssymb}

\begin{document}
\title{Orbital stability of standing waves of some $m$-coupled
nonlinear Schr\"odinger equations}
\date{}
\author{ Hichem Hajaiej}

\maketitle

\begin{abstract}
We extend the notion of orbital stability to systems of nonlinear
Schr\"odinger equations, then we prove this property under suitable
assumptions of the local nonlinearity involved.
\end{abstract}
\section{Intoduction} In [4], the author has studied the following
Cauchy problem :
$$\left\{ \begin{array}{ll}
i\partial_t \Phi_1 + \Delta \Phi_1 +
h_1(x,|\Phi_1|^2,...,|\Phi_\ell|^2) \Phi_1 &=0\\
\vdots \qquad\vdots  \qquad\vdots\\
i\partial_t \Phi_\ell + \Delta \Phi_\ell + h_\ell(x,
|\Phi_1|^2,...,|\Phi_\ell|^2)\Phi_\ell&=0\\
 \Phi_j(0,x) = \Phi^0_j(x) &\mbox{ for } 1 \leq j \leq
 \ell\end{array}\right.
 \eqno{(1.1)}$$
 $\Phi^0_j : \mathbb{R}^N \rightarrow \mathbb{C}, h_j : \mathbb{R}
 \times \mathbb{R}^\ell_+ \rightarrow \mathbb{R}$ and
 $\Phi_j : \mathbb{R} \times \mathbb{R}^N \rightarrow \mathbb{C}$
 (1.1) has numerous applications in engineering and physics. It
 appears in the study of spatial solitons in nonlinear wave guides,
 the theory of Bose-Einstein condensates, optical pulse propagation
 in briefringent fibers, interactions of $m$-wave packets,
 wavelength division multiplexed optical systems, see [3] and
 references therein. Physically $\Phi_j$ is the jth component of the
 beam in Kerr-like photorefractive media. In these contexts it is
 always possible to write (1.1) in a compact vectorial form :
 $$\left\{ \begin{array}{ll}
 i\frac{\partial\vec{\Phi}}{\partial t} &= \hat{E}'(\vec{\Phi})\\
 \vec{\Phi}(0,x) &=\vec{\Phi}^0 = (\Phi^0_1,...,\Phi^0_\ell)
 \end{array}\right.$$
 where $$\hat{E}(\vec{\Phi}) = \frac{1}{2} \left\{ |\nabla\vec{\Phi}|^2_2
 - \int H(x, \Phi_1,..,\Phi_\ell)\right\}\eqno{(1.3)}$$
 $H$ is such that :
 $$\frac{\partial H}{\partial s_j} (x, s_1,...,s_\ell) = 2h_j(x, s^2_1
 ,...,s^2_\ell) s_j
 \eqno{(1.4)}$$
 Thus
 $$\frac{\partial H}{\partial s_j}(x, s_1,...,s_\ell) =
 \frac{\partial H}{\partial s_j}
 (x, |s_1|,...,|s_\ell|)\eqno{(1.5)}$$
 Note that when $\ell = 1$, $ H(x,s) =
 \displaystyle{\int^{s^2}_0}h(x,t)dt$.\\
 A soliton or standing wave of (1.1) is a solution of (1.1) having
 the particular form : $\vec{\Phi}(t,x) =
 (\Phi_1(t,x),...,\Phi_\ell(t,x))$ where $\Phi_j(t,x) = u_j(x)
 e^{-i\lambda_jt}$ ; $\lambda j$ are real numbers.\\
 Hence $\vec{u} = (u_1,...,u_\ell)$ is a solution of the following
 $m \times m$ elliptic eigenvalue problem :
 $$\left\{ \begin{array}{ll}
 \Delta u_1 + h_1(x, u^2_1,...,u^2_\ell)u_1+\lambda_1u_1 &=0\\
 \vdots&\vdots \quad \vdots\\
 \Delta u_\ell + h_\ell(x, u^2_1,...,u^2_\ell) u_\ell + \lambda_\ell
 u_\ell &=0
 \end{array}\right.\eqno{(1.6)}$$
 when $\ell = 1$, (1.6) becomes :
 $$\Delta w + h(x,w^2) w + \lambda w = 0\eqno{(1.7)}$$
 where $w \in H^1( \mathbb{R}^N,\mathbb{C})$ ; which can be written
 as a $2 \times 2$ real elliptic system for $(u,v)$ where $w = (u,v) = u +
 iv$ ; namely
 $$\left\{ \begin{array}{ll}
 \Delta u + h(x, u^2+v^2) u + \lambda u &= 0\\
 \Delta v + h(x, u^2 +v^2)v + \lambda v &= 0
 \end{array}\right.\eqno{(1.8)}$$
where $v \equiv 0$ ; (1.8) leads to the scalar equation
$$\Delta u + h(x,u^2)u + \lambda u = 0.\eqno{(1.9)}$$
(1.9) constitues in itself an important chapter of nonlinear
analysis in which  many brilliant mathematicians as Berger,
Cazenave, Berestyski, Nehari and Lions, have intensively
contributed.  $A$ special attention was addressed to the  case
$h(x,s^2) = |s|^{p-1}$. The famous concentration-compactness
principle was built up by Lions
to study the orbital stability of standing waves of (1.9) [8],[1].\\
In the scalar setting, there are two approaches to determine the
orbital stability of standing waves of (1.1) . The first one reduces
this question to the checking of the strict inequality
$\displaystyle{\frac{d}{d\lambda} \int} u^2_\lambda < 0$ for certain
solutions $u_\lambda$ of (1.9). For non-autonomous equations, it is
hard to establish conditions on the nonlinearity $h$ ensuring the
latter monotonicity property, [7] and references therein. In the
vectorial setting, it does not seem possible to extend this approach
for (1.6). The second alternative exploits the hamiltonian structure
of (1.1) when $\ell = 1$ via the characterization of standing waves
as constrained minimum. We will adapt this approach to generalize
the notion of orbital stability of standing waves of (1.1). We will
then establish stability of the latter particular solutions under
general assumptions on $H$ including the most relevant physical
situations where $H(x, \vec{s})$ converges to a function
$H^\infty(x,\vec{s})$ that depends periodically on $x$.\\ Before
formulating the notion of orbital stability of standing waves of
(1.1), let us first introduce some useful notation :
\begin{eqnarray*}
\vec{H} &=& H^1( \mathbb{R}^N,\mathbb{C}) \times...\times
H^1(\mathbb{R}^N,\mathbb{C}) ; H^1(\mathbb{R}^N, \mathbb{C}) = H\\
\vec{H}^1 &=& H^1( \mathbb{R}^N,\mathbb{R}) \times ...\times H^1(
\mathbb{R}^N,\mathbb{R}) ; H^1( \mathbb{R}^N,\mathbb{R}) = H^1\\
\mbox{ For } z &=& (u,v) \;; |z|^2_H = |z|^2_2 + |\nabla z|^2_2\\
|z|^2_2 &=& |u|^2_2 + |v|^2_2 \;;\quad |\nabla z|^2_2 = |\nabla
u|^2_2 + |\nabla v|^2_2
\end{eqnarray*}
$|\;|_p$ denotes the usual norm on $L^p( \mathbb{R}^N,\mathbb{R}) =
L^p$
$$\vec{z} = (z_1,...,z_\ell) = ((u_1,v_1),...,(u_\ell,v_\ell)) = (\vec{u},\vec{v})$$
where
$$z_j = u_j + iv_j = (u_j,v_j).$$
The modulus of the vector $\vec{z}$, denoted by $|\vec{z}|$ is the
vector
$$|\vec{z}| = (|\vec{z}_1|,...,|\vec{z}_\ell|) \;;\quad
|z_j| = (u^2_j + v^2_j)^{1/2}$$ Let  us now define the following
functionals : $\hat{E} : \vec{H} \rightarrow \mathbb{R}$ and $E :
\vec{H}^1 \rightarrow \mathbb{\mathbb{R}}$.
\begin{eqnarray*}
\hat{E}(\vec{z}) = \hat{E}(\vec{u},\vec{v}) &=& \frac{1}{2}
\{|\nabla \vec{z}|^2_2 - \int H(x, |\vec{z}|)\}\\
&=& \frac{1}{2} \{\sum^\ell_{j=1} \{|\nabla u_j|^2_2 + |\nabla
v_j|^2_2\} - \int H(x,|\vec{z}|)\}\\
&=& \frac{1}{2} \{\sum^\ell_{j=1} |\nabla u_j|^2_2 + |\nabla
v_j|^2_2 - \int
H(x,(u^2_1+v^2_1)^{1/2},...,(u^2_\ell+v^2_\ell)^{1/2}\}
\end{eqnarray*}
$$E(u) = \vec{E}(u,0) = \frac{1}{2} \{|\nabla \vec{u}|^2_2 - \int H(x,|\vec{u}|)\}$$
For $c_1,..,c_\ell > 0$, we set $c^2 =
\displaystyle{\sum^\ell_{i=1}}c^2_i$ and :
$$\vec{S}_c = \{\vec{z} \in \vec{H} : |z_i|^2_2 = c^2_i\quad 1 \;;\leq i \leq \ell\}$$
$$S_c = \{\vec{u} \in \vec{H}^1 : |u_i|^2_2 = c^2_i\quad 1 \leq i \leq \ell\}$$
$$\hat{I}_{c_1,...,c_\ell} = \inf \{\hat{E}(\vec{z}) : \vec{z} \in \vec{S}_c\}$$
and
$$I_{c_1,...,c_\ell} = \inf\{E(\vec{u}) ; \vec{u} \in S_c\}$$
$$\hat{O}_c = \{\vec{z} \in \hat{S}_c : \hat{E}(\vec{z}) =
\hat{I}_{c_1,...,c_2}\}$$
From now on we fix $c_1,...,c_\ell > 0$ and $c^2 =
\sum^\ell_{i=1}c^2_i$.\\
Following the definition in the scalar setting, we will say that
$\hat{O}_c$ is stable if it is not empty and  :
$$\left\{ \begin{array}{l}
\forall\; \vec{w} \in \hat{O}_c\mbox{ and } \forall\; \varepsilon >
0, \exists\; \delta > 0\mbox{ such that }\\
\mbox{ for any }  \vec{\Phi}_0 \in \vec{H} \mbox{ such that }
|\vec{\Phi}_0 - \vec{z}|_{\vec{H}} < \delta, \mbox{ it follows that
}\\
\inf_{\vec{z} \in \hat{O}_c} |\vec{\Phi}(t,.)-\vec{w}|_{\vec{H}} <
\varepsilon \quad \forall\; t\in \mathbb{R}
\end{array}\right.\eqno{(1.10)}$$
$\vec{\Phi}(t,.)$ designs the solution of (1.1) corresponding to the
initial condition $\vec{\Phi}_0$. Hence  we take advantage of the
recent result established in [4], in which the author has determined
assumptions on $h_j$ ensuring the existence and uniqueness of global
solutions of (1.1). Under slight  modifications of Theorem 2.11 and
Theorem 3.1 of [4], we have the following result.\\
{\bf Theorem 0.1} : Let $H : \mathbb{R} \times \mathbb{R}^\ell_+
\rightarrow \mathbb{R}$ be a Carath\'eodory function such that:\\
$(H_0)$ There exist $K > 0$ and $0 < \ell_1 < \frac{4}{N}$ such that
$$0 \leq H (x, \vec{s}) \leq K(|\vec{s}|^2 + |\vec{s}|^{\ell_1+2})$$
for any $x \in \mathbb{R}^N, \vec{s} \in \mathbb{R}^\ell_+$.\\
$(H_1)$ If $N \geq 2$, there exist constants $c' > 0$ and $\alpha
\in [0, \frac{4}{N-2})$ ; for $N \geq 3, \alpha \in [0,\infty)$ for
$N=2$ such that
$$|h_j(x,|\vec{s}|^2) s_j - h_j (x,|\vec{r}|^2)r_j|\leq c'\{1 +
 |\vec{s}|^\alpha + |\vec{r}|^\alpha\}
|\vec{s}-\vec{r}|$$ for all $1 \leq j \leq \ell, \vec{r},\vec{s} \in
\mathbb{R}^\ell_+$.

If $N = 1$, for any $R > 0$, there exists a constant $L(R) > 0$ such
that $|h_j(x,|\vec{s}|^2)s_j - h_j(x,|\vec{r}|^2) r_j|\leq L(R)
|\vec{s}-\vec{r}|$ for all $\vec{s}, \vec{r} \in \mathbb{R}^\ell_+$
such that $|\vec{r}|+|\vec{s}| \leq R$.\\
Then for every $\vec{\Phi}_0 \in \vec{H}$, the initial value, the
Cauchy problem (1.1) has a unique solution $\vec{\Phi} \in
C(\mathbb{R},\vec{H}) \cap C^1( \mathbb{R}, (\vec{H}^{-1})$.
Furthermore $\displaystyle{\sup_{t\in \mathbb{R}}}
|\vec{\Phi}(t,.)|_{\vec{H}} < \infty$ ; and we have conservation of
 charges and energy ; namely
$$(C_1)\quad |\Phi_j(t,.)|_2 = |\Phi^j_0|_2\quad \forall\; 1\leq j
\leq \ell\mbox{ and }
\forall\; t \in \mathbb{R}$$
$$(C_2)\qquad \hat{E}(\vec{\Phi}(t,.)) = \hat{E}(\vec{\Phi}_0)\quad
\forall\; t \in \mathbb{R}.$$ Once one knows that (1.1) admits a
unique solution, it is worth to argue by contradiction to establish
(1.10) :\\
Suppose that $\hat{O}_c$ is not stable, then either $\hat{O}_c$ is
empty or :\\ There exist $\vec{w} \in \hat{O}_c, \varepsilon_0 > 0$
and a sequence $\{\vec{\Phi}_0^n\} \in \vec{H}$ such that :
$$|\vec{\Phi}^n_0 - \vec{w}|_{\vec{H}} \rightarrow 0\mbox{ as }
 n \rightarrow \infty \mbox{ but }
\inf_{\vec{z}\in \hat{O}_c} |\vec{\Phi}^n(t_n,.) -
\vec{z}|_{\vec{H}} \geq \varepsilon_0\}\eqno{(1.11)}$$ for some
sequence $\{t_n\}\subset \mathbb{R}$, where $\Phi^n(t_n,.)$ is the
solution
of (1.1) corresponding to the initial condition $\vec{\Phi}^n_0$.\\
Let $\vec{w}_n = \vec{\Phi}^n(t_n,.)$ ; since $\vec{w} \in
\hat{S}_c$ and $\hat{E}(\vec{w}) = \hat{I}_{c_1,...,c_\ell}$ it
follows from the continuity of $|\;|_2$ and $\hat{E}$ on $\vec{H}$
(Proposition 2.1) that : $|\Phi^n_{0,j}|_2 \rightarrow c_j$
$\forall\; 1 \leq j \leq \ell$ and $\hat{E}(\vec{w}_n)
=\hat{E}(\Phi^n_0) = \hat{I}_{c_1,...,c_\ell}$. Thus it follows from
Theorem 0.1 that
$$|w_{n,j}|_2 = |\Phi^n_{0,j}|_2 \rightarrow c_j\quad \forall\; 1 \leq j \leq \ell$$
and $$\hat{E}(\vec{w}_n) = \hat{E} (\vec{\Phi}^n_0) \rightarrow
\hat{I}_{c_1,...,c_\ell}.$$
 If $\{\vec{w}_n\}$ admits a subsequence
converging to an element $\vec{w} \in \vec{H}$ $\vec{w} =
(w_1,...,w_\ell)$ then $|w_j|_2 \rightarrow c_j$ and
$\hat{E}(\vec{w}) = \hat{I}_{c_1,...,c_\ell}$ showing that $\vec{w}
\in \hat{O}_c$ but $\displaystyle{\inf_{\vec{z}\in \hat{O}_c}}
|\vec{\Phi}^n(t_n,.)-\vec{z}|_{\vec{H}} \leq
\vec{w}_n-\vec{w}|_{\vec{H}}$ contradicting (1.11) . Hence to show
the orbital stability of $\hat{O}_c$, one has to prove that
$\hat{O}_c$ is not empty and :
$$\left\{ \begin{array}{l}
\mbox{ Every  sequence } \{\vec{w}_n\} \subset \vec{H} \mbox{ such
that } |w_{n,j}|_2 \rightarrow c_j\\
\mbox{ for } 1 \leq j \leq \ell \mbox{ and } \hat{E} (\vec{w}_n)
\rightarrow \hat{I}_{c_1,...,c_\ell}
\end{array}\right.\eqno{(1.12)} $$
is relatively compact in $\vec{H}$.\\
In the following $\{\vec{w}_n\}$ denotes a sequence satisfying
(1.12). Our objective is to prove that $\{\vec{w}_n\}$ admits a
subsequence converging to an element $\vec{w} \in \vec{H}$. Our line
of attack consists in the following steps :\\
Step 1 : If $\{\vec{w}_n\}$ satisfies (1.12) then the sequence
$$|\vec{w}_n| = (|\vec{w}_{n,1}|,...,|\vec{w}_{n,\ell}|)\quad
 \mbox{ is such that } E(|\vec{w}_n|) \rightarrow I_{c_1,...,c_\ell} \mbox{ and }
 |w_{n,j}|^2_2 \rightarrow c_j$$ In [5], the author has established
assumptions on $H$ ensuring that such a sequence is relatively
compact in $\vec{H}^1$. It can be easily deduced from Theorem 1.1 of
[5] that :
\\
 {\bf Theorem 0.2} : Suppose that $H$ satisfies $(H_0)$, $(H_1)$
 and\\
$(H_2)$ there exists $B > 0$ such that
$$|\partial_j H(x,\vec{s})| \leq B(|\vec{s}| + |\vec{s}|^{\ell_1+1})\quad
\mbox{ for all } x \in \mathbb{R}^N$$ and $\vec{s} \in
\mathbb{R}^\ell_+ ; 1 \leq j \leq \ell$\\
$(H_3)\quad \exists\; \Delta > 0, R > 0, s > 0,
\alpha_1,...,\alpha_\ell > 0, \quad t \in [0,2)$\\
such that $H(x,\vec{s}) > \Delta |x|^{-t}
|s_1|^{\alpha_1}...|s_\ell|^{\alpha_\ell}$ for all $|x| \geq R$ and
$|\vec{s}| < S$ where $N+2 > \frac{N}{2}\alpha + t\;;  \alpha =
\displaystyle{\sum^\ell_{j=1}}\alpha_j$.\\
$(H_4) H(x, \theta_1 s_1,...,\theta_\ell s_\ell) \geq \theta^2_{max}
H(x, s_1,...,s_\ell)$ for all $x \in \mathbb{R}^N$ $s_i \in
\mathbb{R}_ , \theta_i \geq 1$ where $\theta_{max} = \max_{1 \leq j
\leq \ell} \theta_j$.\\
There exists a periodic function $H^\infty(x,\vec{s})$ (i.e,
$\exists\; T \in \mathbb{Z}^N$ such that $H^\infty  (x, +T,\vec{s})
= H^\infty(x,\vec{s}), \forall\; x \in \mathbb{R}^N, \vec{s} \in
\mathbb{R}^\ell_+)$ satisfying $(H_3)$ and such that :\\
$(H_5)$ There exists $0 < \Gamma < \frac{4}{N}$ such that
$$\lim_{|x|\rightarrow \infty} \frac{H(x, \vec{s})-H^\infty(x, \vec{s})}
{|\vec{s}|^2+|\vec{s}|^{\Gamma+2}} = 0\quad \mbox{ uniformly for any
} \vec{s}$$ $(H_6)$ There exist $A', B' > 0$ and $0 < \beta < \ell_1
< \frac{4}{N}$ such that
$$0 \leq H^\infty(x,\vec{s}) \leq
A'(|\vec{s}|^{\beta+2}+|\vec{s}|^{\ell_1+2}$$ and $\forall\; 1 \leq
j \leq \ell$ :
$$\partial_i H^\infty(x,\vec{s}) \leq \vec B'(|\vec{s}|^{\beta+1} +
|\vec{s}|^{\ell_1+1})\quad \forall\; x \in \mathbb{R}^N\;; \vec{s}
\in \mathbb{R}^\ell_+\;.$$ $(H_7)$ There exists $\sigma \in (0,
\frac{4}{N})$ such that :
$$H^\infty(x,\theta_1 s_1,...,\theta_\ell s_\ell) \geq
 \theta^{\sigma+2}_{max} H^\infty (x, s_1,...,s_\ell)$$
for any $\theta_i \geq 1, x \in \mathbb{R}^N, \vec{s} \in
\mathbb{R}^\ell_+$ , where $\theta_{max} = \max_{1 \leq j \leq \ell}
\theta_i$.\\
Then any sequence $\{\vec{u}_n\} \subset \vec{H}^1$ such that
$|u_{n,j}|^2_2 \rightarrow c_j$ and $E(\vec{u}_n) \rightarrow
I_{c_1,...,c_\ell}$ admits a subsequence converging to $\vec{u} \in
S_c$. Using this important information, step 2 consists of :\\
{\bf Step 2} : By the latter, we now know that there exists $\vec{w}
\in \vec{H}^1$ such that $(u^2_{n,j} + v^2_{n,j})^{1/2}$ converges
to
$w_j$ in $H^1$ for any $1 \leq j \leq \ell$.\\
On the other hand, it follows by Proposition 2.2 that $\vec{w}_n =
((u_{n,1},v_{n,1}), (...,(u_{n,\ell}, v_{n,\ell}))$  is bounded in
$\vec{H}$. Hence up to a subsequence, we may suppose that
$$u_{n,j} > u_j\quad \mbox{ and } \quad v_{n,j} > v_j\quad
\forall\; 1 \leq j \leq \ell$$ In this step, we will prove that $w_j
= (u^2_j + v^2_j)^{1/2}$ $\forall\; 1 \leq j \leq \ell$\\
{\bf Step 3} : We will establish some estimates on $|\nabla
\vec{w}_n|^2_2 -|\nabla |\vec{w}_n||^2_2$, which  will enable us to
prove that $w_{n,j} \rightarrow w_j$ $\forall\; 1 \leq j \leq \ell$
which concludes the proof and here is our main result :\\
{\bf Theorem 1.1}. Suppose that $(H_0)$ to $(H_7)$ are satisfied
then for any $c_1,...,c_\ell > 0$, the orbit $\hat{O}_c$ is stable.
\section{Preliminaries}

Following the proof of Lemma 3.1 of [5], we can easily derive the
following proposition.\\
{\bf Proposition 2.1} : Under the hypothesis $(H_0)$, the
functionals $\hat{E}$ and $E$ are continuous and  have the below
properties
\begin{enumerate}
\item There exists a constant $C > 0$ such that
$$\hat{E} (\vec{z}) \geq \frac{1}{4} |\nabla \vec{z}|^2_2 - C(c^2+c^\gamma)$$
for all $\vec{z} \in \hat{S}_c$ and all $c_1,...,c_\ell > 0$ where
$$\gamma = \frac{2(2\ell_1+4-N\ell_1)}{4-N\ell_1} > 2$$
\item For all $c_1,...,c_\ell > 0$, $I_{c_1,...,c_\ell} \geq
\hat{I}_{c_1,...,c_\ell} > - \infty$ and any minimizing sequences
for $I_{c_1,...,c_\ell}$ and $\hat{I}_{c_1,...,c_\ell}$ are bounded
in $\vec{H}^1$ (resp. $\vec{H}$).
\item $(c_1,...,c_\ell) \rightarrow I_{c_1,...,c_\ell}$ is
continuous on $(0, \infty)^\ell$. \\Now for the convenience of the
under, let us recall a classical.
\end{enumerate}
 {\bf Proposition 2.2.}

 Let $u, v \in H^1$, then $(u^2+v^2)^{1/2} \in H^1$ and for $1 \leq i \leq N$
 $$\partial_i (u^2+v^2)^{1/2} = \left\{ \begin{array}{ll}
\frac{u \partial_i u + v\partial_i v}{(u^2+v^2)^{1/2}} &\mbox{ if }
u^2 + v^2 \neq 0\\
0 &\mbox{ otherwise }
 \end{array}\right.$$
 {\bf Proof} For $\varepsilon > 0$, set \begin{eqnarray*}
 \psi^\varepsilon : \mathbb{R}^2 &\longrightarrow& \mathbb{R}\\
 (s_1,s_2) &\longmapsto& (s^2_1 + s^2_2 + \varepsilon^2)^{1/2} -
 \varepsilon
 \end{eqnarray*}
 Clearly $\psi^\varepsilon \in C^1(\mathbb{R}^2,\mathbb{R}),
  \psi^\varepsilon (0,0) = 0$
 and sup $|\nabla \psi^\varepsilon| < \infty_s$,  it then follows by
 [9] that
 $$\int\left\{ (u^2+v^2+\varepsilon)^{1/2} - \varepsilon\right\}
 \partial_i\xi = - \int\;
 \frac{u \partial_i u+v\partial_iv}{(u^2+v^2+\varepsilon^2)^{1/2}} \xi.$$
 for any $\xi \in C^\infty_0$.\\
 Since $0 \leq (u^2 + v^2 + \varepsilon^2)^{1/2}-\varepsilon \leq (u^2+v^2)^{1/2}$
 and
 $$\left|\frac{u\partial_i u+v\partial_iv}{(u^2+v^2+\varepsilon^2)^{1/2}}\right|
 \leq |\partial_i u| + |\partial_i v|,$$
 we obtain
 $$\int \{u^2(x) + v^2(x)\}^{1/2} \partial_i \xi(x) =
  \int \lim_{\varepsilon \rightarrow 0^+} \frac{u(x)
   \partial_i(x)+ v(x)\partial_i v(x)}
  {u^2(x)+v^2(x)+ \varepsilon^2)^{1/2}}\xi(x)$$
  thanks to the dominated convergence theorem.
  \section{Proof of Theorem 1.1}

  Let $\vec{w}_n = (w_{n,1},...,w_{n,\ell}) = (\vec{u}_n,\vec{v}_n)
  = ((u_{n,1},v_{n,1})...,(u_{n,\ell}, v_{n,\ell}))$
  be a sequence in $\vec{H}$ such that $|w_{n,j}|_2 \rightarrow c_j$
  and $\hat{E}(\vec{w}_n) \rightarrow \hat{I}_{c_1,...,c_\ell}$. We
  will prove that $\{\vec{w}_n\}$ has subsequence converging in
  $\vec{H}$\\
  Setting $|\vec{w}_n| = (c_\ell|w_{n,1}|,...,|w_{n,\ell}|)$, it follows
  by Proposition 2.2 that $|w_{n,j}| = (u^2_{n,j}+v^2_{n,j})^{1/2} \in
  H^1$ and for any $1 \leq j \leq \ell$ and $1 \leq i \leq N$
  $$\partial_i |w_{n,j}| = \left\{\begin{array}{ll}
  \frac{u_{n,j}\partial_i u_{n,j}+v_{n,j}\partial_iv_{n,j}}
  {(u^2_{n,j}+v^2_{n,j})^{1/2}} &\mbox{ if
  } u^2_{n,j} + v^2_{n,j} > 0\\
  0 &\mbox{ otherwise }
  \end{array}\right.$$
  On the other hand, by proposition 2.1, the sequence
  $\{\vec{w}_n\}$ is bounded in $\vec{H}$, and hence passing to a
  subsequence, there exists $\vec{w} = (w_1,...,w_\ell) =
  ((u_1,v_1),...,(u_\ell,v_\ell)) \in \vec{H}$ such that
  $$\left\{ \begin{array}{l}
  $$\forall\; 1 \leq j \leq \ell\quad u_{n,j} >
  u_j, v_{n,j} > v_j\mbox{ and }\\

  \displaystyle{\lim_{n\rightarrow \infty}\int}|\nabla u_{n,j}|^2 +
  |\nabla v_{n,j}|^2\quad \mbox{ exists }\end{array}\right.\eqno{(3.1)}$$
  Now
 $$ \begin{array}{ll}
 \hat{E}(\vec{w}_n) -E(|\vec{w}_n|) &= \frac{1}{2}
 \{|\nabla \vec{w}_n|^2_2 - |\nabla|\vec{w}_n| |^2_2\}\\
 &= \displaystyle{\frac{1}{2} \sum^\ell_{j=1}}|\nabla w_{n,j}|^2_2
 - |\nabla (u^2_{n,j}+v^2_{n,j})^{1/2}|^2_2\\
 &= \displaystyle{\frac{1}{2} \sum^\ell_{j=1} \sum^N_{i=1}
  \frac{(u_{n,j}\partial_i v_{n,j}
 - v_{n,j} \partial_i u_{n,j})^2}{u^2_{n,j}+v^2_{n,j}}} \geq 0\qquad (3.2)
  \end{array}$$
  Proving that
  $$\hat{I}_{c_1,...c_\ell} = \lim_{n\rightarrow \infty}\hat{E}(\vec{w}_n)\geq
  \limsup E(|\vec{w}_n|)\eqno{(3.3)}$$
  But
  $$|w_{n,j}|^2_2 = | |w_{n,j}||^2_2 = c^2_{n,j} \rightarrow c^2_j
  \eqno{(3.4)}$$
  It follows by the continuity property of $I_{c_1,...,c_\ell}$ proved
  in Proposition 2.1 that we have :
  $$\lim \hat{E} (\vec{w}_n) \geq \liminf I_{c_{n,1},...,c_{n,\ell}} =
  I_{c_1,...,c_\ell} \geq \hat{I}_{c_1,...,c_\ell}.$$
  Hence
  $$\lim_{n\rightarrow + \infty} \hat{E} (\vec{w}_n) = \lim_{n\rightarrow \infty}
  E(|\vec{w}_n|) = I_{c_1,...,c_\ell} = \hat{I}_{c_1,...,c_\ell}\eqno{(3.5)}$$
  (3.2) and (3.5) imply that
  $$\forall\; 1 \leq j \leq \ell \;\;\lim_{n\rightarrow \infty}
  \int|\nabla u_{n,j}|^2 + |\nabla v_{n,j}|^2 -
  |\nabla (u^2_{n,j}+v^2_{n,j})^{1/2} = 0\eqno{(3.6)}$$
  Thus it follows form (3.1) that :
  $$\lim_{n\rightarrow \infty} \int|\nabla u_{n,j}|^2 +
   |\nabla v_{n,j}|^2 = \lim_{n\rightarrow + \infty} \int
   |\nabla (u^2_{n,j} + v^2_{n,j})^{1/2}|^2
   \eqno{(3.7)}$$
   which is equivalent to say that :
   $$\lim_{n\rightarrow \infty} |\nabla \vec{w}_n|^2_2 =
   \lim_{n\rightarrow \infty}|\nabla|\vec{w}_n||^2_2\eqno{(3.8)}$$
   (3.4) and (3.5) imply using Theorem 0.2  that
   $|\vec{w}_n|$ is relatively compact in $\vec{H}^1$. Thus there
   exists $w_j \in H^1$ such that
   $$\left\{ \begin{array}{l}
   (u^2_{n,j} + v^2_{n,j})^{1/2} \mbox{ converges to } w_j \mbox{ in
   } H^1 \mbox{ and }\\
   |w_j|_2 = c_j\quad \forall\; 1 \leq j\leq \ell
   \end{array}\right.\eqno{(3.9)}$$
   and $E(w_1,...,w_\ell) = I_{c_1,...,c_\ell}$.\\
   Our purpose now is to prove that $w_j = (u^2_j + v^2_j)^{1/2}$
   ($u_j$ and $v_j$ are given in (3.1)).\\
   Using (3.1), it follows that $u_{n,j} \rightarrow u_j$ and $v_{n,j} \rightarrow
   v_j$ in $L^2(B(0,R))$. Furthermore a straigthfoward computation
   enables us to prove that.
   $$[(u^2_{n,j}+ v^2_{n,j})^{1/2} - (u^2_j + v^2_j)^{1/2}]^2 \leq |u_{n,j}-u_j|^2 + |v^2_{n,j}-v_j|^2$$
   from which we deduce that :
   $$(u^2_{n,j} + v^2_{n,j})^{1/2} \longrightarrow (u^2_j+v^2_j)^{1/2} \mbox{ in }
    L^2(B(0,R))$$
    for all $R > 0$. But $(u^2_{n,j} + v^2_{n,j})^{1/2} \rightarrow
    w_j$ in $L^2$, thus we certainly have that $(u^2_j+v^2_j)^{1/2} =
    w_j$ $\forall\; 1 \leq j \leq \ell$.\\
    On the other hand $|w_{n,j}|_2 = ||w_{n,j}||_2 \rightarrow c_j =
    |w_j|_2$, hence we are done if we prove that $$\lim_{n\rightarrow \infty}
    |\nabla w_{n,j}|^2_2 \rightarrow |\nabla w_j|^2_2 \quad
    \mbox{ for any } 1 \leq j \leq \ell.$$
    Form (3.7) we have that $\displaystyle{\lim_{n\rightarrow \infty}} |\nabla w_{n,j}|^2_2 =
    \displaystyle{\lim_{n\rightarrow \infty}} |\nabla
    |w_{n,j}||^2_2$ and
    $$\lim_{n\rightarrow \infty} |\nabla |w_{n,j}||^2_2 = |\nabla |w_j||^2_2.$$
    Hence
    $$|\nabla w_j|^2_2 \leq \lim |\nabla |w_{n,j}||^2_2
     = |\nabla |w_j||^2_2 .\eqno{(3.9)}$$
     Finally replacing $w_{n,j}$ by $w_j$ in (3.2), we see that
     $$|\nabla w_j|^2_2 \geq |\nabla |w_j||^2_2 \quad
     \forall\; 1 \leq j \leq \ell .\eqno{(3.10)}$$
     By (3.1), we know that $w_{n,j} \rightarrow w_j$ in $H$, thus
     $w_{n,j} > w_j$ ; which completes the proof of
     Theorem 1.1.\vspace{5mm}\\
 {\bf References}

 [1] T. Cazenave, P.L. Lions, Orbital stability of standing waves
 for some nonlinear Schrödinger equations, Comm. Math. Phys. 85, p 549 - 561 (1982).

[2] T. Cazenave, An introduction to nonlinear Schrödinger equations,
Textos de Metodos Matematicos, Rio de Janeiro, 1996, third edition.

[3] H. Hajaiej, Symmetric ground states solutions of m-coupled
nonlinear Schrödinger equations, Nonlinear Analysis: Methods, Theory
and Applications, Vol 71, 2, (2009).

[4] H. Hajaiej, On Schrödinger Systems with Local and Nonlocal
Nonlinearities - Part I, under review.

[5] H.Hajaiej, Existence of Minimizers of a class of
multi-constrained variationnal problems in the absence of
compactness, symmetry and monotonicity, under review.

[6] H.Hajaiej, C.A.Stuart, On the variational approach to the
stability of standing waves for the nonlinear Schrödinger equation,
Adv Nonlinear Studies, 4 (2004), 469-501.

[7] J.B. McLeod, C.A. Stuart, W.C. Troy, Stability of standing waves
for some nonlinear Schr¨odinger equations, J. Diff. Int. Equats., 16
(2003), 1025-1035.

[8]  P.L. Lions, The concentration - compactness principle in the
calculus of variations. The locally compact case: Part 1, p 109 -
145, Part 2, p 223 - 283, Ann. Inst. H.Poincar´e, Vol 1, n 4, 1984.

[9] Jost J¨urgen, Postmodern analysis, Univertex, Springer, 1998.
\end{document}